\documentclass[11pt,reqno]{amsart}
\usepackage{amsmath,amsthm,amsfonts,amssymb,amscd,amstext}
\usepackage[ansinew]{inputenc}
\usepackage[dvips]{graphicx}
\usepackage{psfrag,euscript}

\newcommand{\qand}{\quad\text{and}\quad}

\theoremstyle{plain}
\newtheorem{maintheorem}{Theorem}
\newtheorem{maincorollary}[maintheorem]{Corollary}

\newtheorem{theorem}{Theorem}[section]

\newtheorem{corollary}[theorem]{Corollary}
\newtheorem{lemma}[theorem]{Lemma}

\theoremstyle{definition}

\newtheorem{definition}{Definition}

\newcommand{\RR}{{\mathbb R}}

\newcommand{\TT}{{\mathbb T}}
\newcommand{\DD}{{\mathbb D}}

\newcommand{\vfi}{\varphi}

\newcommand{\real}{{\mathbb R}}

\renewcommand{\epsilon}{\varepsilon}
\newcommand{\dist}{\operatorname{dist}}

\newcommand{\closure}{\operatorname{clos}}

\newcommand{\F}{\EuScript{F}}

\newcommand{\re}{{\mathbb R}}
\newcommand{\V}{\EuScript{V}}
\newcommand{\U}{\EuScript{U}}
\newcommand{\G}{\EuScript{G}}
\newcommand{\X}{\EuScript{X}}

\title
{On the volume of singular-hyperbolic sets}

\author{J. F. Alves, V. Ara\'ujo, M. J. Pacifico, V. Pinheiro}

\date{\today}
\begin{thanks} {J.F.A. and V.A. were partially supported by
    FCT through CMUP, POCI/MAT/61237/2004, and grants of
    FCT(Portugal). M.J.P. was partially supported by
    CNPq-Brazil/FAPERJ-Brazil/Pronex Dynamical Systems.
    V.P.  was partially supported by PADCT/CNPq(Brazil).}
\end{thanks}

\begin{document}

\address{José F. Alves, Centro de Matemática da
  Universidade do Porto
Rua do Campo Alegre 687, 4169-007 Porto, Portugal}
\email{jfalves@fc.up.pt} 

\address{Vítor Araújo, Centro de Matem\'atica da
  Universidade do Porto, Rua do Campo Alegre 687, 4169-007
  Porto, Portugal,
  \emph{and also} Instituto de Matem\'a\-tica,
Universidade Federal do Rio de Janeiro,
C. P. 68.530, 21.945-970
  Rio de Janeiro, RJ-Brazil}
\email{vitor.araujo@im.ufrj.br \textrm{and} vdaraujo@fc.up.pt}

\address{Maria Jos\'e Pacifico,
Instituto de Matem\'atica,
Universidade Federal do Rio de Janeiro,
C. P. 68.530, 21.945-970 Rio de Janeiro, Brazil}
\email{pacifico@im.ufrj.br}

\address{Vilton  Pinheiro, Departamento de Matem\'atica, Universidade Federal da Bahia\\
Av. Ademar de Barros s/n, 40170-110 Salvador, Brazil.}
\email{viltonj@ufba.br}


\renewcommand{\subjclassname}{\textup{2000} Mathematics Subject Classification}


\begin{abstract}
  An attractor $\Lambda$ for a $3$-vector field $X$ is
  singular-hyperbolic if all its singularities are
  hyperbolic and it is partially hyperbolic with volume
  expanding central direction.  We prove that $C^{1+\alpha}$
  singular-hyperbolic attractors, for any $\alpha>0$, always
  have zero volume, extending an analogous result for
  uniformly hyperbolic attractors. The same result holds for
  a class of higher dimensional singular attractors.
  Moreover, we prove that if $\Lambda$ is a
  singular-hyperbolic attractor for $X$ then either it has
  zero volume or $X$ is an Anosov flow. We also present
  examples of $C^1$ singular-hyperbolic attractors with
  positive volume. In addition, we show that $C^1$
  generically we have volume zero for $C^1$ robust classes
  of singular-hyperbolic attractors.
\end{abstract}

\maketitle

\setcounter{tocdepth}{1}

 \tableofcontents

\section{Introduction}

The uniform hyperbolic theory of Dynamical Systems was introduced in
the 60's by Smale~\cite{Sm67} and nowadays the understanding of
uniformly hyperbolic systems is fairly complete. However, open sets
of systems, such as perturbations of the Lorenz flow \cite{Lo63} or
perturbations of the geometric Lorenz flows \cite{ABS77,GW79} fail
to be uniformly hyperbolic.  These flows are not uniformly
hyperbolic because they present equilibria accumulated in a robust
way by regular orbits.  The study of the dynamical properties of
this kind of flows has been the object of many works in the field
and lead to the concept of a \emph{singular-hyperbolic set}, which
is a compact invariant set $\Lambda$ for a $3$-flow such that  all
its singularities are hyperbolic and it is partially hyperbolic with
volume expanding central direction \cite{MPP04,MPP98}; a precise
definition of all these concepts will be given below.

A\emph{ singular-hyperbolic attractor} is an attracting set which is singular-hyperbolic and contains 
a dense orbit. 
The first examples of singular hyperbolic sets included the Lorenz
attractor \cite{Lo63,Tu99} and its geometric models
\cite{Gu76,ABS77,GW79,Wil79}, and the singular-horseshoe
\cite{LP86}, besides the uniformly hyperbolic sets themselves.
Many other examples have been found recently, including attractors
arising from certain resonant double homoclinic loops \cite{MPsM}
or from certain singular cycles \cite{MPu97}, and certain models
across the boundary of uniform hyperbolicity \cite{MPP1}.

The next natural step is to understand the dynamical
consequences of singular hyperbolicity.  From the
topological point of view this has been very successful
\cite{MPP99,CMP00,MP01,MP03,MPP04,Mo04,BM04,LMP05}.

Our goal here is to complement the previous topological
results from the measure theoretical point of view. Early
work in this direction was obtained in \cite{Ki88}, where it
was shown that the geometric Lorenz model from \cite{GW79}
has an SRB measure which is stochastically stable.  More
recently, the existence of SRB measures for $C^2$
singular-hyperbolic attractors was proved in
\cite{clm2002,APPV}.  Furthermore, it was shown in
\cite{APPV} that this measure has a disintegration into
conditional measures along the central direction that are
absolutely continuous with respect to the Lesbegue measure,
and moreover the support of the SRB measure is the whole
attractor.

Here we show that singular-hyperbolic attractors of
$C^{1+\alpha}$ flows, for any $\alpha>0$, always have zero
volume, extending both an analogous result for uniformly
hyperbolic attractors \cite{BR75} and a generalization
for partially hyperbolic sets from \cite{alpi2005}. We point out
that since the divergence of both the geometric Lorenz
attractor and the flow of the Lorenz equations is negative,
it is true for these particular systems that the attractors
have zero volume~\cite{robinson2004}. From our proofs we
further deduce that the same zero volume property holds for
a class of higher dimensional singular-attractors, as the
so-called multidimensional Lorenz attractors \cite{BPV97}.

We also prove a dichotomy related with the zero volume property: if
$\Lambda$ is a $C^{1+\alpha}$ singular-hyperbolic attractor for a
$3$-dimensional vector field $X$, then either it has zero volume or
$X$ induces a transitive Anosov flow. In addition, we show that
$C^1$ generically we have zero volume for $C^1$ robust classes of
singular-hyperbolic attractors. We also present examples of $C^1$
singular-hyperbolic attractors with positive volume.




\subsection{Partial hyperbolicity and singular-hyperbolicity}
\label{sec:sing-hyperb}

Let $M$ be a compact boundaryless $d$-dimensional manifold, for some
$d\ge3$, and let ${\X}^r(M)$ be the set of $C^r$ vector fields on
$M$, endowed with the $C^r$ topology, for some $r\ge1$. From now on
we fix some smooth Riemannian structure on $M$ and an induced
normalized volume form $m$ that we call Lebesgue measure. We write
also $\dist$ for the induced distance on $M$. Given $X \in
{\X}^1(M)$, we denote by $(X_t)_{t \in \re}$ the flow induced by
$X$.

We also denote by $\X^{1+}(M)$ the set of all $C^1$ vector fields
$X$ whose derivative $DX$ is H\"older continuous with respect to the
given Riemannian norm,  and we say that $X\in\X^{1+}(M)$ is of class
$C^{1+}$. We obviously have
$$\X^1(M)\supset\X^{1+}(M)\supset\X^r(M), \quad\text{for every $r\ge2$}.$$
Let $\Lambda$ be a compact invariant set of $X\in \X^1(M)$.  We say
that $\Lambda$ is \emph{isolated} if there exists an open set
$U\supset \Lambda$ such that
$$
\Lambda =\bigcap_{t\in \re}X_t(U).
$$
If $U$ above can be chosen such that $X_t(U)\subset U$ for every
$t>0$, then we say that $\Lambda$ is an \emph{attracting
  set}.
We say that an attracting set $\Lambda$ is \emph{transitive} if it
coincides with the $\omega$-limit set of a regular orbit.  An
\emph{attractor} is a transitive attracting set, and a
\emph{repeller} is an attractor for the reversed vector field $-X$.
An attractor, or repeller, is \emph{proper} if it does not coincide
with  the whole manifold.  An invariant set of $X$ is
\emph{non-trivial} if it is neither a periodic orbit nor a
singularity.

\begin{definition}
\label{d.dominado} Let $\Lambda$ be a compact invariant set for $X
\in \X^1(M)$. Given $0 < \lambda < 1$, we say that $\Lambda$ has a
{\em dominated splitting} if the tangent bundle over $\Lambda$ can
be written as a continuous $DX_t$-invariant sum of sub-bundles
\[
T_\Lambda M=E\oplus F,
\]
such that, for some choice of a Riemannian metric $\|\quad\|$, we
have for every $t
> 0$ and every $x \in \Lambda$
\begin{equation}\label{eq.domination}
\|DX_t \mid E_x\| \cdot \|DX_{-t} \mid F_{X_t(x)}\| <  \lambda^t.
\end{equation}
We say that $\Lambda$ is \emph{partially hyperbolic} if it has a
dominated splitting  as in \eqref{eq.domination} 
for which $E$ is \emph{uniformly contracting}, i.e.  for every $t >
0$ and every $x \in \Lambda$ we have \( \|DX_t \mid E_x\| <
\lambda^t. \)
\end{definition}

For $x\in \Lambda$ and $t\in\real$ we let $J_t(x)$ be the absolute
value of the determinant of the linear map
\[
DX_t \mid F_x:F_x\to F_{X_t(x)}.
\]
We say that the sub-bundle $F$ of the partially hyperbolic invariant
set $\Lambda$ is \emph{volume
  expanding} if $J_t(x)\geq \, e^{\lambda t}$ for every
$x\in \Lambda$ and every $t\geq 0$.



\begin{definition}
\label{d.singularset}
Let $\Lambda$ be a compact invariant set for $X \in
\X^r(M)$, $r\ge1$.  We say that $\Lambda$ is a
\emph{singular-hyperbolic set} for $X$ if all the
singularities of $\Lambda$ are hyperbolic, and $\Lambda$ is
partially hyperbolic with volume expanding central direction
either for $X$ or for the reverse flow $-X$.
\end{definition}

We emphasize that \emph{a compact invariant set $\Lambda$ without
singularities for $X$ on a $3$-manifold which is partially
hyperbolic with volume expanding central direction}, either for $X$
or for $-X$  is a \emph{uniformly hyperbolic set of saddle-type}, by
\cite[Proposition 1.8]{MPP04}. In the $3$-dimensional setting this
means that $E$ is one-dimensional and $F$ can be written as a sum of
two invariant one-dimensional sub-bundles $F=[X]\oplus G$, where
$[X]$ is the flow direction and $G$ is uniformly expanding.
 In
general, a compact invariant set $\Lambda$ for $X$ is uniformly
hyperbolic of saddle-type if $$T_\Lambda M=E\oplus [X] \oplus G$$ is
a continuous $DX_t$-invariant splitting with the sub-bundle $E\neq0$
uniformly contracted and the sub-bundle $G\neq0$ uniformly expanded
by $DX_t$ for $ t>0$.

An embedded disk $\gamma\subset M$ is a (local) {\em
  strong-unstable manifold}, or a {\em strong-unstable
  disk}, if $\dist(X_{-t}(x),X_{-t}(y))$ tends to zero
exponentially fast as $t\to+\infty$, for every
$x,y\in\gamma$. Similarly, $\gamma$ is called a (local) {\em
  strong-stable manifold}, or a {\em strong-stable disk}, if
$\dist(X_{t}(x),X_{t}(y))\to0$ exponentially fast as
$n\to+\infty$, for every $x,y\in\gamma$. It is well-known
that every point in a uniformly hyperbolic set possesses a
local strong-stable manifold $W_{loc}^{ss}(x)$ and a local
strong-unstable manifold $W_{loc}^{uu}(x)$ which are disks
tangent to $E_x$ and $G_x$ at $x$ respectively with
topological dimensions $d_E=\dim(E)$ and $d_G=\dim(G)$
respectively. Considering the action of the flow we get the
(global) \emph{strong-stable manifold}
$$W^{ss}(x)=\bigcup_{t>0}
X_{-t}\Big(W^{ss}_{loc}\big(X_t(x)\big)\Big)$$
and the
(global) \emph{strong-unstable manifold}
$$W^{uu}(x)=\bigcup_{t>0}X_{t}\Big(W^{uu}_{loc}\big(X_{-t}(x)\big)\Big)$$
for every point $x$ of a uniformly hyperbolic set. Similar
notions are defined in a straightforward way for
diffeomorphisms. These are immersed submanidfolds with the
same differentiability of the flow or the diffeomorphism.
In the case of a flow we also consider what we call the
\emph{stable manifold} $W^s(x)=\cup_{t\in\real}
X_{t}\big(W^{ss}(x)\big)$ and \emph{unstable
  manifold}
$W^u(x)=\cup_{t\in\real}X_{t}\big(W^{uu}(x)\big)$ for $x$ in
a uniformly hyperbolic set, which are flow invariant.


\subsection{Statement of the results}
\label{sec:statement-results}


Our first result generalizes the results of
Bowen-Ruelle~\cite{BR75} which show that a uniformly
hyperbolic transitive subset of saddle-type for a $C^{1+}$
flow has zero volume.

\begin{maintheorem}
  \label{mthm:C+attractor}
  Let $X\in\X^{1+}(M)$ be a vector field on a 3-dimensional
  manifold $M$.  Then any proper singular-hyperbolic attractor
  or repeller for $X$ has zero volume.
\end{maintheorem}

 Moreover, we obtain the
following dichotomy extending a similar result obtained by two of
the authors in \cite{alpi2005} to the continuous time setting.

\begin{maintheorem}
  \label{mcor:basiczerovolume}
  Let $\Lambda$ be a transitive isolated uniformly
  hyperbolic set of saddle type for $X\in\X^{1+}(M)$ on a
  $d$-dimensional manifold $M$, for some $d\ge3$. Then either
  $\Lambda$ has zero volume or $X$ is a transitive Anosov
  vector field.
\end{maintheorem}

A \emph{transitive Anosov} vector field $X$ is a vector field
without singularities such that the entire manifold $M$ is a
uniformly hyperbolic set of saddle-type. Using these results we can
deduce the following statement which generalizes to
singular-hyperbolic attractors (a class which includes all
transitive uniformly hyperbolic invariant subsets of saddle-type) a
similar one obtained in \cite{alpi2005} for transitive uniformly
hyperbolic sets of $C^{1+}$ diffeomorphisms.

\begin{maintheorem}
  \label{mthm:volposanosov}
  Let $\Lambda$ be a singular-hyperbolic attractor for
  $X\in\X^{1+}(M)$ where $M$ is a $3$-manifold. Then either
  $\Lambda$ has zero volume or $X$ is a transitive Anosov
  vector field.
\end{maintheorem}

%
%

An interesting consequence of the zero volume property of
$C^{1+}$ singular-hyperbolic attractors and uniformly
hyperbolic sets of saddle-type coupled with an extension of
the notion of Axiom A system is as follows.

Recall that a vector field $X\in\X^1(M)$ is \emph{Axiom A} if the
non-wandering set $\Omega(X)$ is both hyperbolic and the closure of
its periodic orbits and singularities. The \emph{spectral
decomposition theorem} (see e.g. \cite{Sh87}) asserts that if $X$ is
Axiom A, then there is a disjoint decomposition
\[
\Omega(X) =  H_1\cup\dots\cup H_n,
\]
where each $H_i$ is a hyperbolic basic set for $X$ for
$i=1,\dots,n$. Following~\cite{MPP99} we say that a vector field
$X\in\X^1(M)$ is \emph{singular Axiom A} if there is a finite
disjoint decomposition
\[
\Omega(X) = \Lambda_1\cup\dots\cup \Lambda_n,
\]
where each $\Lambda_i$ is a hyperbolic basic set, a
singular-hyperbolic attractor or a singular-hyperbolic repeller for
$i=1,\dots,n$.  We remark that singular Axiom~A systems play an
important part in the generic description of $C^1$-flows on
$3$-manifolds as proved in~\cite{MP03}.  The \emph{topological
basin} of an attracting set $\Lambda$ is the set
  \[
  W^s(\Lambda)= \big\{ x\in M : \lim_{t\to+\infty}\dist\big(
  X_t(x) , \Lambda\big) =0 \big\}.
\]
A straightforward consequence of Theorems~\ref{mthm:C+attractor}
and~\ref{mcor:basiczerovolume} and the well-known pairwise disjoint
decomposition
\[
M=W^s(\Lambda_1)\cup\dots\cup W^s(\Lambda_n)
\]
from \cite[Chpt. 2, Lemma 2.2]{Sh87} is the following.

\begin{maincorollary}
  \label{mcor:qtpatrator}
  If $X\in\X^{1+}(M)$ is a singular Axiom A vector field on a
  $3$-manifold $M$, then
\[
\cup\big\{ W^s(\Lambda) : \Lambda \mbox{ is a
  hyperbolic attractor or a singular-hyperbolic attractor} \big\}
\]
has full Lebesgue measure in $M$.
\end{maincorollary}

A construction in \cite{BPV97} shows that there exist $C^1$
open sets of vector fields exhibiting proper robust
attractors containing hyperbolic singularities with any
number $k\ge2$ of expanding directions on $d$-manifolds $M$
with $d=k+3$. Moreover as will be shown in
Section~\ref{sec:robust-attr-high} these attractors satisfy
the conclusions of Theorem~\ref{mthm:C+attractor}.  We say
that an attractor $\Lambda$ for $X\in\X^1(M)$ with isolating
neighborhood $U$ is \emph{robust} if there is a $C^1$
neighborhood $\U$ of $X$ in $\X^1(M)$ such that the maximal
invariant set
\[
\Lambda_Y(U)=\bigcap_{t\in\RR} \overline{Y_t(U)}
\]
is transitive for all $Y\in\U$. In the case of 3-manifolds this
implies that the set $\Lambda_Y(U)$ is a singular-hyperbolic
attractor \cite{MPP04}.

\begin{maincorollary}
  \label{mcor:singattractor-highdim}
  A multidimensional proper singular attractor for
  the class constructed in \cite{BPV97} has volume zero.
\end{maincorollary}

The properties stated above are not valid in general in the
$C^1$ setting: we present in
Section~\ref{sec:posit-volume-sing} an explicit construction
of a geometric Lorenz attractor for a flow of class $C^1$ on
any $3$-manifold with positive Lebesgue measure.

However we are able to extend the previous results to some locally
generic subsets in the $C^1$ topology. Recall that a subset $\G$ of
a set $\U$ is \emph{generic} if it may be written as a countable
intersection of open and dense subsets of~$\U$. Since $\X^1(M)$ is a
Baire space, generic subsets of a given open subset $\U$ of
$\X^1(M)$ are dense in $\U$.

Concrete examples of such open sets on 3-manifolds are given by
robust singular-hyperbolic attractors, see e.g.~\cite{MPu97,MPP04}
which comprise the Lorenz attractor~\cite{Tu99}, the geometric
Lorenz attractors~\cite{ABS77,GW79},  attractors
arising from certain resonant double homoclinic loops \cite{MPsM}
or from certain singular cycles \cite{MPu97}, and certain models
across the boundary of uniform hyperbolicity \cite{MPP1}. 

\begin{maintheorem}
  \label{mcor:3drobust}
  Let $\Lambda$ be a robust attractor for $X\in\X^1(M)$ on a
  $3$-manifold $M$ with isolating neighborhood $U$. Then
  there is a $C^1$-neighborhood $\U$ of $X$ and a
  $C^1$-generic set $\G\subset\U$ such that $\Lambda_Y(U)$
  has volume zero for all $Y\in\G$.
\end{maintheorem}

As mentioned above the $C^1$-open sets of
singular-attractors on $d$-manifolds described
in~\cite{BPV97} are also in the setting of
Theorem~\ref{mthm:C1horseshoe}.

\begin{maintheorem}
  \label{mcor:4drobust}
  Let $\Lambda$ be a robust multidimensional
  singular attractor for a
  vector field $X$  with isolating
  neighborhood $U$ as in~\cite{BPV97}. Then there
  exists a $C^1$-neighborhood $\U$ of $X$ and a
  $C^1$-generic set $\G\subset\U$ such that $\Lambda_Y(U)$ has volume zero for all
  $Y\in\G$.
\end{maintheorem}


This paper is organized as follows. Theorem~\ref{mthm:C+attractor}
is a consequence of Theorem~\ref{mthm:C+partialhorseshoe} proved in
Section~\ref{sec:part-hyperb-zero} together with
Lemma~\ref{le:sdisks-anosov} proved in
Subsection~\ref{sec:posit-volume-sing-1}. Both
Theorems~\ref{mcor:basiczerovolume} and~\ref{mthm:volposanosov} are
proved in Section~\ref{sec:posit-volume-vers}. Then in
Section~\ref{sec:posit-volume-sing} we present an example of a
positive volume $C^1$ singular-hyperbolic attractor and in
Section~\ref{sec:robust-attr-high} we show that certain classes of
multidimensional robust singular-attractors are in the setting of
Theorem~\ref{mthm:C+partialhorseshoe}. Finally we show that
Theorems~\ref{mcor:3drobust} and~\ref{mcor:4drobust} are
consequences of Theorem~\ref{mthm:C1horseshoe} proved in
Section~\ref{sec:zero-volume-c1}.

\medskip
\noindent
\textbf{Acknowledgment:} We thank Marcelo Viana for helpfull
discussions on the statements and proofs contained in this work.


\section{Partial hyperbolicity and zero volume on  $C^{1+}$ flows}
\label{sec:part-hyperb-zero}

The following result plays a crucial role in the proof of
Theorem~\ref{mthm:C+attractor}.

\begin{theorem}
  \label{mthm:C+partialhorseshoe}
  Let $X$ be a $C^{1+}$ flow on a $d$-dimensional manifold
  with $d\ge3$ and $\Lambda$ be a partially hyperbolic
  invariant subset such that
\begin{enumerate}
\item[($\star$)] $\Lambda\cap \gamma$ does not contain $d_E$-disks
  for any strong-stable disk $\gamma$.
\end{enumerate}
Then $\Lambda$ has zero volume.
\end{theorem}

The proof of the above theorem is a consequence of the
following result for compact invariant subsets of $C^{1+}$
diffeomorphisms $f$ with dominated decomposition, whose
proof we present in Subsections~\ref{s.holder}
and~\ref{s.local}.


Before we state the result let us recall the notion of
\emph{dominated splitting for a diffeomorphism} $f$ over a
compact $f$-invariant set $\Lambda$, which is very similar
to the one given by Definition~\ref{d.dominado} since we
need only exchange (\ref{eq.domination}) by
\begin{equation}
  \label{eq:domination-diffeo}
\|Df \mid E_x\| \cdot \|Df^{-1} \mid F_{x}\| < \lambda
\end{equation}
for all $x\in\Lambda$.  Analogously \emph{partial
  hyperbolicity for a diffeomorphism} $f$ is given by a
dominated decomposition $E\oplus F$ over a compact invariant
subset $\Lambda$ with uniform contraction along the
direction $E$.

\begin{theorem}\label{t:disco1}
Let \( f: M\to M \) be a  \( C^{1+} \)
  diffeomorphism and let $\Lambda\subset M$ be a partially
  hyperbolic set with positive volume.  Then $\Lambda$
  contains a strong-stable disk.
\end{theorem}


Let us now prove Theorem~\ref{mthm:C+partialhorseshoe}
  using Theorem~\ref{t:disco1}.
Let $\Lambda$ be a partially
  hyperbolic compact invariant set for a flow
  $X\in\X^{1+}(M)$ where $M$ is a $d$-manifold with $d\ge3$.
  Assume that condition ($\star$) is satisfied by $\Lambda$.

Arguing by contradiction, if $m(\Lambda)>0$ then setting $f=X_{1}$,
the time-one diffeomorphism induced by the vector field $X$, we see
that $\Lambda$ is in the setting of Theorem~\ref{t:disco1}.

Hence there exists some strong-stable disk $\gamma$ for $f$
contained in $\Lambda$ with dimension $d_E$, which is a
strong-stable disk for the flow $X_t, t>0$ and contradicts
property ($\star$).  This contradiction shows that $m(\Lambda)=0$
and proves Theorem~\ref{mthm:C+partialhorseshoe}.


\subsection{Pre-balls and bounded
  distortion}\label{s.holder}

Here we give some preparatory results for the proof of
Theorem~\ref{t:disco1}. We fix continuous extensions of the two
bundles $E$ and $F$ to some neighborhood $U$ of $\Lambda$, that we
denote by $\tilde{E}$ and $\tilde{F}$. We do not require these
extensions to be invariant under $Df$. Given $0<a<1$, we define the
{\em center-unstable cone field
  $\left(C_a^{F}(x)\right)_{x\in U}$ of width $a$\/} by
\begin{equation}
\label{e.cucone} C_a^{F}(x)=\big\{v_1+v_2 \in \tilde{E}_x\oplus
\tilde{F}_x \mbox{\ such\ that\ } \|v_1\| \le a\cdot \|v_2\|\big\}.
\end{equation}
We define the {\em stable cone field
  $\left(C_a^{E}(x)\right)_{x\in U}$ of width $a$\/} in a
similar way, just reversing the roles of the bundles in
(\ref{e.cucone}). We fix $a>0$ and $U$ small enough so that,
up to slightly increasing $\lambda<1$, the domination
condition~(\ref{eq:domination-diffeo}) remains valid for any
pair of vectors in the two cone fields:
\begin{equation}\label{domina2}
    \|Df(x)u\|\cdot\|Df^{-1}(f(x))v\|
\le\lambda\cdot\|u\|\cdot\|v\|
\end{equation}
for every $u\in C_a^{E}(x)$, $v\in C_a^{F}(f(x))$, and any
point $x\in U\cap f^{-1}(U)$. Note that the unstable cone field is
positively invariant:
\[
Df(x) C_a^{F}(x)\subset C_a^{F}(f(x)),
\]
whenever $x,f(x)\in U$. Indeed, the domination (\ref{domina2})
together with the invariance of $F=\tilde{F} \mid \Lambda$
imply that
\[
Df(x) C_a^{F}(x) \subset C_{\lambda a}^{F}(f(x))
                  \subset C_a^{F}(f(x)),
\]
for every $x\in \Lambda$. This extends to any $x\in U\cap
f^{-1}(U)$ just by continuity. Analogously the stable cone field is
negatively invariant:
\[
Df^{-1}(x) C_a^{E}(x)\subset C_a^{E}(f^{-1}(x)),
\]
whenever $x,f(x)\in U$.

If $a>0$ is taken sufficiently small in the definition of
the cone fields, and we choose $\delta_1>0$ also small so
that the $\delta_1$-neighborhood of $\Lambda$ should be
contained in~$U$, then by continuity
\begin{equation}
\label{delta1}
\|Df(y) u \| \le \lambda^{-1/2}\cdot
\big\|Df\mid E_{x}\big\| \cdot\|u\|,
\end{equation}
whenever $x\in \Lambda$, $\dist(x,y)\le \delta_1$, and $u\in
C^{E}_a(y)$.

We say that an embedded $C^1$ submanifold $N\subset U$ is {\em
tangent to the stable cone field\/} if the tangent subspace to $N$
at each point $x\in N$ is contained in $C_a^{E}(x)$.  Then, by the
domination property~\eqref{eq:domination-diffeo}, $f^{-1}(N)$ is
also tangent to the stable cone field, if it is contained in $U$. In
particular, if $N,f^{-1}(N),\dots, f^{-k}(N)\subset U$, then $Df^{k}
\mid T_{f^{-k}(x)N}$ is a $\lambda^{k/2}$-contraction
by~(\ref{delta1}), since $\big\|Df\mid E_{x}\big\|<1$ by partial
hyperbolicity. Thus, denoting by $\dist_N$ the \emph{distance along
$N$} given by the length of the shortest smooth curve connecting two
given points inside $N$, we obtain

\begin{lemma} \label{l.contraction} Let $\Delta\subset U$ be
  a $C^1$ disk of radius $\delta<\delta_1$ tangent to the
  stable cone field.
There exists $n_0\ge1$ such that for $n\ge n_0$ and $x\in \Delta$
with $\dist_\Delta(x,\partial \Delta)\ge \delta/2$ there is a
neighborhood $V_n$ of $x$ in $\Delta$ such that $f^{-n}$ maps $V_n$
diffeomorphically onto a disk  of radius $\delta_1$ around
$f^{-n}(x)$. Moreover,
\[
\dist_{f^{-n+k}(V_n)}(f^{-n+k}(y),f^{-n+k}(z)) \le
\lambda^{k/2}\cdot
\dist_{f^{-n}(V_n)}(f^{-n}(y),f^{-n}(z))
\]
for every $1\le k \le n$ and every $y, z\in V_n$.
\end{lemma}

We shall sometimes refer to the sets \( V_n \) as
\emph{pre-balls}.  The next corollary is a consequence of
the contraction given by the previous lemma, together with
some H\"older control of the tangent direction which can be
found in \cite[Corollary 2.4, Proposition 2.8]{ABV00}

\begin{corollary} \label{p.distortion} There exists $C>1$ such that given
$\Delta$ as in Lemma~\ref{l.contraction}
 and given any pre-ball $V_n\subset
\Delta$ with $n\ge n_0$, then for all $y,z\in V_n$
$$
\frac{1}{C} \le \frac{|\det Df^{-n} \mid T_y \Delta|}
                     {|\det Df^{-n} \mid T_z \Delta|}
            \le C.
$$
 \end{corollary}


\subsection{A local unstable disk inside $\Lambda$}\label{s.local}

%

Assume that $\Lambda$ has positive volume and given an
embedded disk $\Delta$ in $M$ denote by $m_\Delta$ the
measure naturally induced by the volume form $m$ on
$\Delta$.  Choosing a $m$ density point of $\Lambda$, we
laminate a neighborhood of that point into disks tangent to
the stable cone field. Since the relative Lebesgue measure
of the intersections of these disks with $\Lambda$ cannot be
all equal to zero, we obtain some disk $\Delta$ intersecting
$\Lambda$ in a positive $m_\Delta$ subset. Hence, in the
setting of Theorem~\ref{t:disco1}, we assure that there is a
disk~$\Delta$ tangent to the stable cone field
intersecting $\Lambda$ in a positive $m_\Delta$ subset.
%
%
%
%
Let $H=\Delta\cap\Lambda$.

\begin{lemma} \label{l.discao} There exist an infinite sequence of integers
$1\le k_1<k_2<\cdots$ and, for each $n\in\mathbb{N}$,  a disk
$\Delta_n$ of radius $\delta_1/4$ tangent to the stable cone field
such that the relative Lebesgue measure of $f^{-k_n}(H)$ in
$\Delta_n$ converges to~1 as $n\to\infty$.
\end{lemma}

\begin{proof} Let $\epsilon>0$ be some small number.
  Let $K$ be a compact subset of $H$ and $A$ be an open
  neighborhood of $H$ in $\Delta$ such that
\[
m_\Delta(A\setminus K)<\epsilon{m_\Delta(K)}.
\]
Choose $n$ sufficiently large so that for each $x\in K$ we
have $V_x\subset A$, where $V_x$ is the pre-ball associated
to $n$. This pre-ball is mapped diffeomorphically by $f^{-n}$
onto a ball $B_{\delta_1}(f^{-n}(x))$ of radius $\delta_1$
around $f^{-n}(x)$ tangent to the stable cone field. Let
$W_{x}\subset V_x$ be the pre-image of the ball
$B_{\delta_1/4}(f^{-n}(x))$ of radius $\delta_1/4$ under this
diffeomorphism. By compactness we have
\[
K\subset
W_{x_1}\cup...\cup W_{x_m},
\]
for some $x_1,...,x_m\in K$.
Let $I$ be a maximal set of $\{1,...,m\}$ such that for
$i,j\in I$ with $i\neq j$ we have $W_{x_i}\cap
W_{x_j}=\emptyset$. By maximality, each $W_{x_j}$, $1\le
j\le m$, intersects some $W_{x_i}$ with $i\in\mathcal U$.
Hence $\{V_{x_i}\}_{i\in I}$ is a covering of~$K$. By
Corollary~\ref{p.distortion} there is a uniform constant
$\theta>0$ such that
\[
\frac{m_\Delta(W_{x_i})}{m_\Delta(V_{x_i})}\ge\theta,
\quad\mbox{for every }i\in I.
\]
Hence
\[
m_\Delta\big(\cup_{i\in I}W_{x_i}\big) = \sum_{i\in I}
m_\Delta(W_{x_i}) \ge
    \sum_{i\in I}\theta  m_\Delta(V_{x_i})
\ge \theta  m_\Delta(K).
\]
Setting
\[
\rho=\min\left\{\frac{m_\Delta(W_{x_i}\setminus
K)}{m_\Delta(W_{x_i})}\colon i\in I\right\},
\]
we have
\begin{eqnarray*}
    \varepsilon   m_\Delta(K)
&\ge&
m_\Delta(A\setminus K) \\
&=& m_\Delta\big(\cup_{i\in I}W_{x_i}\setminus
K\big) \\
     &=& \rho  m_\Delta\big(\cup_{i\in I}W_{x_i}\big) \\
    &=& \rho  \theta m_\Delta(K).
\end{eqnarray*}
This implies that $\rho<\varepsilon/\theta$. Since
$\varepsilon>0$ can be taken arbitrarily small, increasing
$n$
%
we may take $W_{x_i}$ such that the relative Lebesgue
measure of $K$ in $W_{x_i}$ is arbitrarily close to 1. Then,
by the bounded distortion provided by
Corollary~\ref{p.distortion}, the relative Lebesgue measure
of $f^{-n}(H)\supset f^{-n}(K)$ in $f^{-n}(W_{x_i})$, which
is a disk of radius $\delta_1/4$ around $f^{-n}(x_i)$
tangent to unstable cone field, can be made arbitrarily
close to $1$.
\end{proof}

Let us now prove that there is a strong-stable disk
of radius $\delta_1/4$ inside~$\Lambda$. Let $(\Delta_n)_n$ be the
sequence of disks given by Lemma~\ref{l.discao}, and consider
$(x_n)_n$ the sequence of points at which these disks are centered.
Up to taking subsequences, we may assume that the centers of the
disks converge to some point $x$. By Ascoli-Arzela, these disks
converge to some disk $\Delta_\infty$ centered at $x$. By
construction, every point in $\Delta_\infty$ is accumulated by some
iterate of a point in $H\subset\Lambda$, and so
$\Delta_\infty\subset \Lambda$.

Note that each $\Delta_n$ is contained in the $k_n$-iterate
of $\Delta$, which is a disk tangent to the stable cone
field. The domination property implies that the angle
between $\Delta_n$ and $E$ goes to zero as $n\to\infty$,
uniformly on $\Lambda$. In particular, $\Delta_\infty$ is
tangent to $E$ at every point in
$\Delta_\infty\subset\Lambda$. By Lemma~\ref{l.contraction},
given any $k\ge 1$, then $f^{k}$ is a
$\sigma^{k/2}$-contraction on $\Delta_n$ for every large
$n$. Passing to the limit, we get that $f^{k}$ is a
$\sigma^{k/2}$-contraction on $\Delta_\infty$ for every
$k\ge1$. In particular, we have shown that the subspace
$E_x$ is uniformly contracting for $Df$. The fact that
$T_\Lambda M=E\oplus F$ is a dominated splitting implies
that any contraction $Df$ may exhibit along the
complementary direction $F_x$ is weaker than the contraction
in the $E_x$ direction.  Then, by \cite{Pe76}, there exists
a unique strong-stable manifold $W^{ss}_{loc}(x)$ tangent to
$E$ and which is contracted by the positive iterates of
$f$. 
Since $\Delta_\infty$ is contracted by every~$f^{k}$, and  all its
positive iterates are tangent to stable cone field, then
$\Delta_\infty$ is contained in $W^{ss}_{loc}(x)$.

This completes the proof of Theorem~\ref{t:disco1}.


\section{Positive volume versus transitive Anosov flows}
\label{sec:posit-volume-vers}

In this section we prove
Theorems~\ref{mthm:C+attractor},
\ref{mcor:basiczerovolume}
and~\ref{mthm:volposanosov}.

\subsection{Positive volume transitive hyperbolic sets and Anosov flows}
\label{sec:posit-volume-trans}

We start by proving Theorem~\ref{mcor:basiczerovolume}. Let
$\Lambda$ be a transitive uniformly hyperbolic set for
$X\in\X^{1+}(M)$ such that $m(\Lambda)>0$, where $M$ is a
$d$-manifold, for some $d\ge3$. 

\begin{lemma}
  \label{le:unifhiposvol-ssdentro}
If there exists a point $x\in\Lambda$ in the interior of
$W^{ss}_{loc}(x)\cap\Lambda$, then $\Lambda\supset W^{ss}(y)$ for
all $y\in\Lambda$. Moreover, the set $W^u(\Lambda)$ formed by the
union of all unstable manifolds through points of $\Lambda$ is an
open neighborhood of $\Lambda$.
\end{lemma}

Here the interior of $W^{ss}_{loc}(x)\cap\Lambda$ is taken with
respect to the topology of the disk $W^{ss}_{loc}(x)$.
The proof follows \cite[Lemma 2.16]{MPP04} and \cite[Lemma
2.8]{APPV} closely.

\begin{proof}
  Let $x\in\Lambda$ be such that
  $x$ is in the interior of $W^{ss}(x)\cap\Lambda$. Let
  $\alpha(x)\subset\Lambda$ be its $\alpha$-limit set. Then
\begin{equation}\label{eq.ssinside}
W^{ss}(z)\subset\Lambda \quad\text{for every } z\in \alpha(x),
\end{equation}
since any compact part of the strong-stable manifold of $z$
is accumulated by backward iterates of any small neighborhood
of $x$ inside $W^{ss}(x)$. Here we are using that the
contraction along the strong-stable manifold, which becomes
an expansion for negative time, is uniform.

Clearly the invariant set $\alpha(x)\subset\Lambda$ is
uniformly hyperbolic.

It also follows from \eqref{eq.ssinside} that the union
\[
S=\bigcup_{y\in\alpha(x)} W^{ss}(y) \quad\mbox{or}\quad
S=W^{ss}(\alpha(x))
\]
of the strong-stable manifolds through the points of $\alpha(x)$ is
contained in $\Lambda$. By continuity of the strong-stable manifolds
and the fact that $\alpha(x)$ is a closed set, we get that $S$ is
also closed.  Again $S$ is a uniformly hyperbolic~set.

We claim that $W^{u}(S)$, the union of the unstable manifolds of the
points of $S$,  is an open set. To prove this, we note that $S$
contains the whole stable manifold $W^{s}(z)$ of every $z\in S$:
this is because $S$ is invariant and contains the strong-stable
manifold of~$z$. Now the union of the strong-unstable manifolds
through the points of $W^{s}(z)$ contains a neighborhood of $z$.
This proves that $W^u(S)$ is a neighborhood of $S$. Thus the
backward orbit of any point in $W^u(S)$ must enter the interior of
$W^u(S)$. Since the interior is, clearly, an invariant set, this
proves that $W^u(S)$ is open, as claimed.

Finally, consider any backward dense orbit in $\Lambda$ of a point
that we call~$w$.  On the one hand $\alpha(w)=\Lambda$.  On the
other hand, $X_{-t}(w)$ must belong to $W^u(S)$ for some $t>0$, and
so $\alpha(w)\subset S$ by invariance.  This implies that
$\Lambda\subset S$ and since $S\subset\Lambda$ by construction, we
see that $\Lambda=S$.
\end{proof}

\begin{proof}[Proof of Theorem~\ref{mcor:basiczerovolume}]
Assume that $m(\Lambda)>0$. Then Theorem~\ref{t:disco1}
applied to the map $f=X_1$ and to the set $\Lambda$ with dominated
decomposition given by the splitting
\[
 E\oplus \big( [X] \oplus  F\big),
\]
ensures that there exists a strong-stable disk $\gamma$ contained in
$\Lambda$. Analogously applying Theorem~\ref{t:disco1} to $f=X_{-1}$
and to the set $\Lambda$ with dominated decomposition given by the
splitting
\[
\big( E\oplus [X] \big) \oplus  F,
\]
we get a strong-unstable disk $\delta$ contained in $\Lambda$.

Now the existence of $\gamma$ enables us to use
Lemma~\ref{le:unifhiposvol-ssdentro} and deduce that $\Lambda$
contains the strong-stable manifolds of each of its points and that
$W^u(\Lambda)$ is an open neighborhood of $\Lambda$. In the same
way, using Lemma~\ref{le:unifhiposvol-ssdentro} for the flow
generated by $-X$, from the existence of $\delta$ we deduce that
$\Lambda$ contains the strong-unstable manifolds of all of its
points, that is $W^u(\Lambda)\subset\Lambda$.

But since $W^u(\Lambda)$ is an open neighborhood
of $\Lambda$, we conclude that $\Lambda$ is simultaneously
open and closed in $M$. Hence $\Lambda=M$ by
connectedness. This shows that the whole of $M$ is a
transitive uniformly hyperbolic set for $X$ and completes
the proof of Theorem~\ref{mcor:basiczerovolume}.
\end{proof}


\subsection{Positive volume singular-hyperbolic sets and
  Anosov flows}
\label{sec:posit-volume-sing-1}

Now we prove Theorems~\ref{mthm:C+attractor}
and~\ref{mthm:volposanosov}.  For that we need some
preliminary results which show in particular that transitive
singular-hyperbolic sets satisfy condition ($\star$).

In what follows $X$ is a vector field in $\X^{1+}(M)$ and
$M$ is a $3$-manifold.

\begin{lemma}
\label{le:sdisks-unif-hyp}
Let $\Lambda$ be a transitive partially hyperbolic invariant
set for $X$ with volume expanding central direction. Then
\begin{itemize}
\item either $W^{ss}(x)\cap\Lambda$ contains no
  strong stable disks for all $x\in\Lambda$,
\item or $\Lambda$ is a uniformly hyperbolic set (and in
  particular $\Lambda$ does not contain singularities).
\end{itemize}
\end{lemma}

\begin{proof}
  Let us suppose that there exists $x\in\Lambda$ such that
  $x$ is in the interior of $W^{ss}(x)\cap\Lambda$. Let
  $\alpha(x)\subset\Lambda$ be its $\alpha$-limit set. Then
  we have (\ref{eq.ssinside}).
It follows that $\alpha(x)$ does not contain any singularity.
Indeed, \cite[Theorem~B]{MPP04} proves that the strong-stable
manifold of each singularity meets $\Lambda$ only at the
singularity. Therefore by \cite[Proposition~1.8]{MPP04} the
invariant set $\alpha(x)\subset\Lambda$ is uniformly hyperbolic.

As in the proof of Lemma~\ref{le:unifhiposvol-ssdentro} we
have that
\[
S=W^{ss}(\alpha(x))\subset\Lambda
\]
and that $S$ is closed. Again we see that $S$ does not contain any
singularity~$\sigma$, for otherwise we would have
$W^{ss}(\sigma)\supset W^{ss}(z)$ for some $z\in\alpha(x)$ which
by~(\ref{eq.ssinside}) would contradict \cite[Theorem~B]{MPP04}.
Thus $S$ is a uniformly hyperbolic~set.

Then $W^{u}(S)$ is also an open set as in the proof of
Lemma~\ref{le:unifhiposvol-ssdentro}. Since we are assumming
that $\Lambda$ is transitive, again by the same arguments in
the proof of Lemma~\ref{le:unifhiposvol-ssdentro} we get
that $\Lambda=S$. This shows that $\Lambda$ is uniformly
hyperbolic and in particular it does not contain
any singularity of $X$. 
\end{proof}

\begin{proof}[Proof of Theorem~\ref{mthm:C+attractor}]
  Note that if $\Lambda$ is a transitive singular-hyperbolic
  set for $X$, then since $\Lambda$ contains singularities,
  by Lemma~\ref{le:sdisks-unif-hyp} we have that $\Lambda$
  satisfies property ($\star$) in the statement of
  Theorem~\ref{mthm:C+partialhorseshoe}.  Hence $\Lambda$
  has zero volume, thus concluding the proof of
  Theorem~\ref{mthm:C+attractor} since a singular-hyperbolic
  attractor is transitive by definition.
\end{proof}

We can obtain a stronger conclusion if we further assume
that $\Lambda$ is an attractor.

\begin{lemma}
  \label{le:sdisks-anosov}
  Let $\Lambda$ be a singular-hyperbolic attractor for $X$.
  Then
\begin{itemize}
\item either $W^{ss}(x)\cap\Lambda$ contains no
  strong stable disks for all $x\in\Lambda$,
\item or $\Lambda=M$ is a uniformly hyperbolic set and $X$
  is a transitive Anosov vector field.
\end{itemize}
\end{lemma}

\begin{proof}
  Assume that there exists $x\in\Lambda$ such that
  $x$ is in the interior of $W^{ss}(x)\cap\Lambda$.
  From Lemma~\ref{le:sdisks-unif-hyp} we know that there
  exists a uniformly hyperbolic set $S\subset\Lambda$ such
  that $W^u(S)$ is an open neighborhood of $S$. Moreover we
  also have that $\Lambda=S$.

However if $\Lambda$ is an attractor, then
$W^u(S)\subset\Lambda$ and so we get $\Lambda=W^u(S)$. Hence
$\Lambda$ is closed and also open. The connectedness of $M$
implies that $\Lambda=S=M$. In particular $X$ has no
singularities and the whole of $M$ admits a uniformly
hyperbolic structure with a dense orbit, thus $X$ is a
transitive Anosov vector field.
\end{proof}


\begin{proof}[Proof of Theorem~\ref{mthm:volposanosov}]
  Let $\Lambda$ be a singular-hyperbolic attractor for a
  $C^{1+}$ vector field $X$ on a $3$-manifold. If
  $m(\Lambda)>0$ then according to Theorem~\ref{t:disco1} we
  get that there exists some strong-stable disk
  $\gamma$ contained in $\Lambda$.

Hence since $\Lambda$ does not satisfy the first alternative of
Lemma~\ref{le:sdisks-anosov} we conclude that $\Lambda=M$
and so $X$ is an Anosov vector field. This concludes the
proof of Theorem~\ref{mthm:volposanosov}.
\end{proof}


\section{A positive volume singular-hyperbolic attractor}
\label{sec:posit-volume-sing}

To construct an example of a 3-dimensional $C^1$ flow exhibiting a
singular-hyperbolic attractor with positive volume we start with the
construction of a \emph{dynamically
  defined Cantor set with positive one-dimensional Lebesgue
  measure}.

This construction is carried out in detail in \cite[Section
4.2]{PT93} giving (after a trivial change of coordinates) a
$C^1$ 2-to-1 surjective map
$$\vfi:[-1/2,a]\cup[b,1/2]\to[-1/2,1/2],$$ such that both
$\vfi\mid[-1/2,a]$ and $\vfi\mid[b,1/2]$ are diffeomorphisms
onto $[-1/2,1/2]$, and $-1/2<a<0<b<1/2$ are fixed. The map
$\vfi$ is such that $\vfi'$ is \emph{not of bounded
  variation nor H\"older continuous}.  Moreover, the
construction is performed in such a way that the maximal positive
invariant set (dynamically defined Cantor set)
\[
K=\bigcap_{n\ge0} \vfi^{-n}([-1/2,1/2]) \quad\mbox{satisfies}\quad
\lambda(K)>0,
\]
where we denote by $\lambda$ the standard Lebesgue measure
on the real line.

\begin{figure}[htbp]
\begin{center}
  \includegraphics[height=5cm]{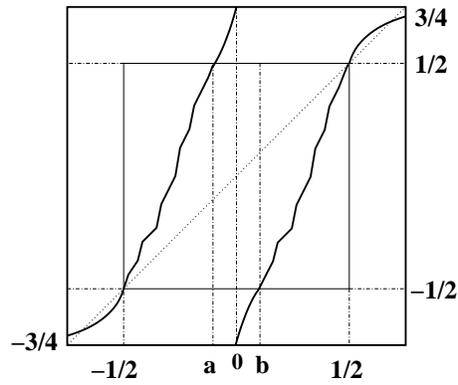}
\end{center}
  \caption{The one-dimensional Lorenz-like map.}
  \label{fig:lorenz-like}
\end{figure}

Now we adapt this map $\vfi$ so that it becomes a \emph{Lorenz-like}
map, which can be used to define a geometric Lorenz flow in the
sense of \cite{ABS77,GW79}. We may assume without loss that both
$\vfi\mid[-1/2,a]$ and $\vfi\mid[b,1/2]$ are increasing. We now
consider a  $C^1$ extension $\phi$ of $\vfi$ to
$[-3/4,3/4]\setminus\{0\}$ satisfying (see
Figure~\ref{fig:lorenz-like}):
\begin{enumerate}
\item $\phi\mid[-3/4,0)$ and $\phi\mid(0,3/4]$ are
  increasing;
\item $-3/4<\phi(-3/4)$ and $\phi(3/4)<3/4$;
\item $\lim_{x\to0^-}\phi(x)=3/4$ and $\lim_{x\to0^+}\phi(x)=-3/4$;
\item $\lim_{x\to0^-}\phi'(x)=+\infty$ and
  $\lim_{x\to0^+}\phi'(x)=-\infty$.
\end{enumerate}
Note also that the maximal positive invariant subset for
$\phi$ in $[-3/4,3/4]$ is the same set $K$ as before and so has
positive one-dimensional Lebesgue measure.

\begin{figure}[htbp]
\begin{center}
\psfrag{S}{\Large$\Sigma$}
  \includegraphics[height=5cm,width=8cm]{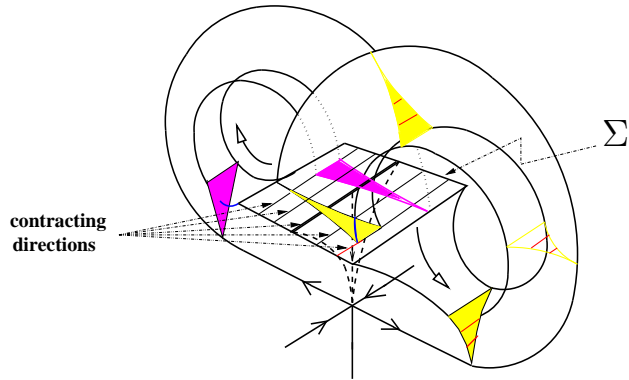}
\end{center}
\caption{The geometric Lorenz attractor with the contracting
  directions on the cross-section $\Sigma$.}
  \label{fig:geom-lorenz}
\end{figure}

Now this map can be used as a basis for the standard
construction of a flow $X$ exhibiting a \emph{geometric
  Lorenz attractor} as explicitly described in
\cite[Chapter 7, Section 3.2]{robinson2004} and sketched in
Figure~\ref{fig:geom-lorenz}.
If we denote by $\Lambda$ the attractor obtained for the vector
field $X$ just constructed, then the set $K$ can now be interpreted
as a projection of the set $\Lambda\cap\Sigma$ through the stable
leaves crossing the cross-section $\Sigma$ (drawn in
Figure~\ref{fig:geom-lorenz}).

Since the projection along this stable leaves is smooth,
because it coincides with the usual projection along lines
with constant coordinate in Euclidean space, this shows that
the two-dimensional Lebesgue measure (area) of
$\Lambda\cap\Sigma$ must be positive.

In what follows, given $A\subset M$ and $[a,b]\subset \re$ we shall
write
$$X_{[a,b]}(A)=\{X_t(x)\,\colon\, a\leq t \leq b \qand x\in A\}.$$

Finally using that $\Sigma$ is a cross-section for the flow,
we consider a small flow box through $\Sigma$ and the set
$\Lambda_\Sigma^\epsilon=
X_{[-\epsilon,\epsilon]}(\Lambda\cap\Sigma)$ for small
$\epsilon>0$ is contained in $\Lambda$ by flow invariance
and has positive volume, since it can be seen as a direct
product of $[-\epsilon,\epsilon]\times(\Lambda\cap\Sigma)$,
because $X_t$ is a diffeomorphism for all $t\in\RR$.

Therefore we conclude that $m(\Lambda)\ge
m\big(\Lambda_\Sigma^\epsilon\big) >0$.


\section{Robust attractors in higher dimensions}
\label{sec:robust-attr-high}

Here we briefly describe the construction of singular-attractors in
\cite{BPV97} with any number of expanding dimensions, and show that
this class of attractors satisfies condition ($\star$).

Consider a ``solenoid'' constructed over a uniformly expanding map
$f:\TT^k\to\TT^k$ of the $k$-dimensional torus, for some $k\ge2$.
That is, let $\DD$ be the unit disk on $\RR^2$ and consider a smooth
embedding $F:\TT^k\times\DD\to\TT^k\times\DD$ of $N=\TT^k\times\DD$
into itself, which preserves and contracts the foliation
$$\F^s=\big\{\{z\}\times\DD: z\in\TT^k\big\},$$ and moreover the natural
projection $\pi:N\to\TT^k$ on the first factor conjugates $F$ to
$f$: $\pi\circ F=f\circ\pi$.

\begin{figure}[htbp]
\centering
\includegraphics[width=7cm]{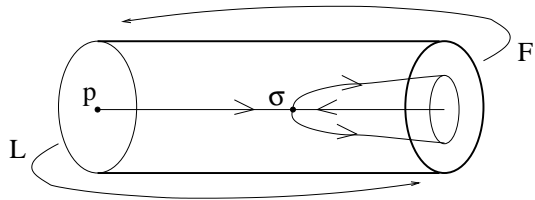}
\caption{A sketch of the construction of a
  singular-attractor in higher dimensions}
\label{fig:sing-attractor-4d}
\end{figure}

Now consider the linear flow over defined on
$M=N\times[0,1]$ given by the vector field $X=(0,1)$ on
$TN\times\RR$. Modify the flow on  a cylinder
$U\times\DD\times[0,1]$ around the orbit of a point
$p=(z,0)\in N$, where $U$ is a neighborhood of $z$ in
$\TT^k$, in such a way as to create a hyperbolic singularity
$\sigma$ of saddle-type with $k$-expanding and $3$
contracting eigenvalues, as depicted in
Figure~\ref{fig:sing-attractor-4d}.

This modified flow defines a transition map $L$ from
$\Sigma_0=\TT^k\times\{0\}$ to $\Sigma_1=\TT^k\times\{1\}$
which through the identification
given by $(1,w)\sim_F (0,F(w))$ defines the return map to
the global cross-section $\Sigma_0$ of a flow $Y$ on the
space $M^F=M/\sim_F$.

In~\cite{BPV97} it is shown that if the expanding rate of $f$ is
sufficiently big, then the set
\[
\Lambda= \bigcup_{T>0} \overline{ \bigcap_{t>T} Y_t(\Sigma_0) }
\]
is a robust partially hyperbolic attractor with
singularities.

To see that $\Lambda$ satisfies condition ($\star$), note first that
the intersection $\Lambda\cap\Sigma_0$ is contained in the
hyperbolic solenoid $\Lambda_0=\cap_{n\ge0} F^n(N)$. Moreover, the
strong-stable manifold of every given point $(z,w,t)\in\Lambda$
(where $(z,w,t)\in\TT^k\times\DD\times[0,1]$) contains the disk
$\{z\}\times\DD\times\{t\}$ which is a leaf of $\F^s$. Hence
$$W^{s}_{loc}(z,w,0)\cap\Lambda\subset \big(\{z\}\times\DD\cap
\Lambda_0\big)\times\{0\},$$
 and this last set is a \emph{Cantor
set}, since it is the intersection of a strong-stable manifold with
a uniformly hyperbolic solenoid. Thus
$W^{s}_{loc}(z,w,0)\cap\Lambda$ does not contain any $2$-disk.

Since $\Sigma_0$ is a global cross-section, any other strong-stable
leaf $W$ is such that $W\cap\Lambda$ does not contain any $2$-disk
$D$, for otherwise by the flow invariance of $\Lambda$ the first
return to $\Sigma_0$ of the points of $D$ would be contained in a
strong-stable leaf and would contain an open set, contradicting the
previous paragraph. This shows that this class of robust
multidimensional singular attractors is in the setting of
Theorem~\ref{mthm:C+partialhorseshoe} and so they have zero volume
in the $C^{1+}$ setting.


\section{Zero volume in the $C^1$ generic setting}
\label{sec:zero-volume-c1}

Here we prove Theorems~\ref{mcor:3drobust}
and~\ref{mcor:4drobust} as a consequence of the following
result.

\begin{theorem}
  \label{mthm:C1horseshoe}
  Let $\Lambda$ be an isolated partially hyperbolic set
  satisfying condition ($\star$)  for $X\in\X^1(M)$ on a
  $d$-dimensional manifold $M$ with $d\ge3$. Given an isolating
  neighborhood $U$ of $\Lambda$, let $\U\subset\X^1(M)$ be  such that
\( \Lambda_Y(U)\) is partially hyperbolic and also satisfies
condition ($\star$) for all $Y\in\U$.

If $\U$ is $C^1$-open, then there exists a generic set $\G\subset\U$
such that $\Lambda_Y(U)$ has volume zero  for all $Y\in\G$.
\end{theorem}

\begin{proof}
  Let $\Lambda$ be an isolated partially hyperbolic
  invariant compact subset for a flow $X\in\X^r(M)$,
for some  $r\ge1$, such that $\Lambda$ satisfies condition
  ($\star$). We always write $U$ for the isolating
  neighborhood of $\Lambda$.

  We consider the sets
\begin{itemize}
\item $\U=\{Y\in\X^1(M): \Lambda_Y(U) \mbox{ is partially
    hyperbolic satisfying ($\star$)}\}$ which we assume is a
  $C^1$ open subset of $\X^1(M)$;
\item $\V=\{ Y\in\X^2(M): \Lambda_Y(U) \mbox{ is partially
    hyperbolic satisfying ($\star$)}\}$;
\item $\U_\epsilon=\{ Y\in\U : m\big(\Lambda_Y(U)\big)
  <\epsilon\}$.
\end{itemize}
Since every $C^1$ flow $X$ is arbitrarily close to some
$C^2$ flow $Y$ in the $C^1$ topology (see e.g. \cite{PM82})
and we are assumming that $\U$ is $C^1$-open, we conclude
that $\V$ is dense in $\U$ in the $C^1$ topology.

We claim that $\U_\epsilon$ is open and dense in $\U$ in the
$C^1$ topology. After proving this claim the proof of
Theorem~\ref{mthm:C1horseshoe} finishes by setting
$\G=\cap_{n\ge1} \U_{1/n}$. In what follows we prove this
claim.

Let $Y\in\U_\epsilon$ be given. Then for every fixed $T>0$
we set
\[
\Lambda^T_Y=\cap_{t=-T}^T \closure\big(Y_t(U)\big)\subseteq
\Lambda_Y(U) \quad\mbox{and}\quad
\epsilon_1=\epsilon-m\big(\Lambda^T_Y\big)>0.
\]
There exists $\delta>0$ such that
\[
m\Big( B\big( \Lambda^T_Y , \delta \big) \setminus
\Lambda^T_Y \Big) < \frac{\epsilon_1}2.
\]
Let $B^{C_1}(Y,\zeta)$ denote the $C^1$-neighborhood of radius
$\zeta$ around $Y$. Using that $T$ is finite, by continuity and
compactness we find $\zeta>0$ such that
\[
Z\in B^{C_1}(Y,\zeta)\quad\Rightarrow\quad \Lambda_Z^T\subset B\big(
\Lambda^T_Y , \delta \big).
\]
We have that for all $Z\in\U\cap B^{C_1}(Y,\zeta)$
\[
m\big( \Lambda_Z^T \big) \le m \Big( B\big( \Lambda^T_Y ,
\delta \big) \Big) \le
m\big(\Lambda^T_Y\big)+\frac{\epsilon_1}2 \le
\frac{m\big(\Lambda^T_Y\big)+\epsilon}2 <\epsilon.
\]
Since $\Lambda_Z(U)\subset \Lambda_Z^T$, we conclude $m\big(
\Lambda_Z(U) \big) <\epsilon$, for all $z\in\U\cap
B^{C_1}(Y,\zeta)$.

This proves that $\U_\epsilon$ is $C^1$-open. To prove that
$\U_\epsilon$ is $C^1$-dense in $\U$, just observe that
$\V\cap\U\subset\U_\epsilon$ by
Theorem~\ref{mthm:C+partialhorseshoe}. Since $\V$ is
$C^1$-dense in $\U$ this concludes the proof of the claim
and ends the proof of Theorem~\ref{mthm:C1horseshoe}.
\end{proof}

\begin{proof}[Proof of Theorem~\ref{mcor:3drobust}]
  This is a straightforward consequence of Theorem~\ref{mthm:C1horseshoe}
 since a $C^1$ robust attractor for $3$-flows is a
  singular-hyperbolic attractor as shown in \cite{MPP04}.
\end{proof}

\begin{proof}[Proof of Theorem~\ref{mcor:4drobust}]
This is also an immediate consequence of  Theorem~\ref{mthm:C1horseshoe}
since the class of multidimensional singular-attractors
constructed in~\cite{BPV97} is in the setting of
Theorem~\ref{mthm:C+partialhorseshoe}, as shown in
Section~\ref{sec:robust-attr-high}.
\end{proof}



\def\cprime{$'$}

\end{document}